\documentstyle{amsart}

\def\Hom{{\rm Hom}}
\def\Ext{{\rm Ext}}
\def\Proj{{\rm Proj}}
\def\SHom{{\cal H}\!{\it om}}
\def\SExt{{\cal E}\!{\it xt}}
\def\pd{{\rm pd}}

\def\cG{{\cal G}}
\def\cF{{\cal F}}
\def\cN{{\cal N}}
\def\cM{{\cal M}}
\def\cO{{\cal O}}
\def\cE{{\cal E}}
\def\cL{{\cal L}}

\def\mm{{\frak m}}
\def\nn{{\frak n}}

\def\PP{{\Bbb P}}
\def\ZZ{{\Bbb Z}}
\def\NN{{\Bbb N}}
\newsymbol\Bbbk 207C
\def\kk{{\Bbbk}}
\def\HH{{\bf H}}
\def\H{{\rm H}}

\def\ol{\overline}
\def\ul{\underline}
\def\ds{\displaystyle}

\title[Computing $\Ext^{m}(X;\cM,\cN)$]{Computing Global Extension Modules \\ for Coherent Sheaves on a Projective Scheme}
\author[G.~G.~Smith]{{Gregory~G.~Smith}}
\address{Department of Mathematics, University of California, Berkeley, California, 94720-3840 U.S.A.}
\email{{\tt ggsmith@@math.berkeley.edu}}

\newtheorem{lemma}{Lemma}[section]  
\newtheorem{theorem}{Theorem}
\newtheorem{corollary}[lemma]{Corollary}
\newtheorem{proposition}[lemma]{Proposition}

\theoremstyle{definition}
\newtheorem{example}[lemma]{Example}
\newtheorem{algorithm}[lemma]{Algorithm}
\newtheorem{remark}[lemma]{Remark}
\newtheorem*{proof}{Proof}
\newtheorem*{freeres}{Free Resolutions}
\newtheorem*{geometry}{Algebraic Geometry}
\newtheorem*{local}{Local Cohomology}
\newtheorem*{hypercohomology}{Hypercohomology}
\newtheorem*{duality}{Duality Theory}
\newtheorem*{cohomology}{Cohomology of Locally Free Sheaves}
\newtheorem*{extensions}{Global Extensions}
\newtheorem*{proofthm}{Proof of Theorem~\ref{main}}

\begin{document}

\begin{abstract}
Let $X$ be a projective scheme; let $\cM$ and $\cN$ be two coherent $\cO_{X}$-modules.  Given an integer $m$, we present an algorithm for computing the global extension module $\Ext^{m}(X;\cM,\cN)$.  In particular, this allows one to calculate the sheaf cohomology $H^{m}(X,\cN)$ and to construct the sheaf corresponding to an element of the module $\Ext^{1}(X;\cM,\cN)$.  This algorithm can be implemented using only the computation of Gr\"{o}bner bases and syzygies, and it has been implemented in the computer algebra system {\tt Macaulay2}.
\end{abstract}

\maketitle

\section*{Introduction}

\noindent Let $X$ be a projective scheme over the field $\kk$; let $\cM$ and $\cN$ be two coherent $\cO_{X}$-modules.  We write $\Ext^{m}(X;\cM,\cN)$ for the $m$-th global extension module.  There are several equivalent definitions for these modules.  Grothendieck \cite{T} and Hartshorne \cite{Har} define them as the right derived functors of $\cN \mapsto \Hom(X;\cM,\cN)$, where $\Hom(X;\cM,\cN)$ denotes the group of $\cO_{X}$-module morphisms.  Griffiths and Harris \cite{G-H} define them via hypercohomology and discuss the correspondence between the equivalence classes of global extensions of $\cM$ by $\cN$, also called Yoneda Ext, and $\Ext^{m}(X;\cM,\cN)$.  The global extension modules play a central role in duality theory --- again see Hartshorne \cite{Har} or Griffiths and Harris \cite{G-H} for more details.  The purpose of this paper is to provide an algorithm for computing $\Ext^{m}(X;\cM,\cN)$.

To achieve this goal, we must determine how to represent the coherent sheaves $\cM$ and $\cN$.  Because we are working over a projective scheme, every coherent sheaf can be represented by a finitely generated graded $R$-module, where $R$ is the homogeneous coordinate ring of $X$.  By embedding $X$ into $\PP^{n}$, we may identify $R$ with a quotient of the polynomial ring $S = \kk[x_{0},\ldots,x_{n}]$.  This representation is advantageous since it allows one to use the algorithms and techniques of computational algebra --- see Vasconcelos \cite{Vas}.  However, there is no canonical choice for the finitely generated module associated to a given coherent sheaf.  In particular, the natural candidate $\bigoplus_{v \in \ZZ} \Gamma\big(X,\cN(v)\big)$ may fail to be finitely generated.  With this in mind, our problem is more precisely stated as follows:  given two finitely generated graded $R$-modules $M$ and $N$, calculate the $R$-module $\Ext^{m}(X;\cM,\cN)$ where $\cM$ and $\cN$ are the coherent sheaves associated to $M$ and $N$ respectively.

The answer for $m=0$ motivates our general result.  In this case, $\Ext^{0}(X;\cM,\cN)$, which equals $\Hom(X;\cM,\cN)$, can be identified with $\big( \Hom_{R}(M_{\geq r}, N)\big)_{0}$ for $r \gg 0$, that is the set of degree zero $R$-linear maps from a truncation $M_{\geq r}$ to $N$, for sufficiently large $r$.  For a simple example in which truncation arises, we consider $M = R$ and $N = \mm^{d}$ with $d >0$; the second module is the $d$-th power of the irrelevant ideal in $R$.  Both of these modules correspond to the structure sheaf $\cO_{X}$.  Since a morphism from $\cO_{X}$ is determined by the image of the unit global section and $\Gamma(X,\cO_{X}) \cong \kk$, it follows that $\Hom(X;\cO_{X},\cO_{X}) \cong \kk$.  Analogously, when $r \geq d$ we have $\big(\Hom_{R}(R_{\geq r},\mm^{d})\big)_{0} \cong \kk$.  However, for $r < d$ there are no degree zero maps from $R_{\geq r}$ to $\mm^{d}$.  Hence, to evaluate $\Ext^{0}(X;\cM,\cN)$ using this bijection, we must determine an appropriate truncation of $M$.  To provide an algorithm, we need an effective bound $r_{0}$ such that $\Ext^{0}(X;\cM,\cN) \cong \big( \Hom_{R}(M_{\geq r},N)\big)_{0}$ for all $r \geq r_{0}$.

Our solution to the general question involves relating $\Ext^{m}(X;\cM,\cN)$ to the $R$-module $\Ext_{R}^{m}(M_{\geq r},N)$; the second Ext module is calculated in the category of graded $R$-modules.  This method has the added virtue that it allows one to understand $\bigoplus_{v \in \ZZ} \Ext^{m}\big(X;\cM,\cN(v)\big)$ as an $R$-module.  To express our main result, we write ${_{S}N}$ for the $S$-module obtained from $N$ by restriction of scalars, $\ol{a}_{i}({_{S}N})$ for the maximal degree of the minimal generators of the $i$-th syzygy module of ${_{S}N}$, and $\pd_{S}({_{S}N})$ for the projective dimension of ${_{S}N}$.

\renewcommand{\thetheorem}{1}
\begin{theorem} \label{main} 
Let $M$ and $N$ be two finitely generated graded $R$-modules, let $m$ be an integer and set $\ell = \min\{\dim(N),m\}$.  If $r$ is an integer satisfying 
$$ r \geq \max\big\{ \ol{a}_{i}({_{S}N})-i : n-\ell \leq i \leq \pd_{S}({_{S}N}) \big\}-m+1,$$
then we have an isomorphism of graded $R$-modules
$$ \bigoplus_{v \geq 0} \Ext^{m}\big(X;\cM,\cN(v)\big) \cong  \Big(\Ext_{R}^{m}(M_{\geq r},N)\Big)_{\geq 0}.$$
\end{theorem}

\noindent In particular, when $\pd_{S}({_{S}N})< n-\ell$, we have $\max\{\emptyset\}=-\infty$ and the inequality for $r$ is vacuously satisfied.  Furthermore observe that, by replacing $N$ with $N(e)$ for some integer $e$, one can compute the module $\bigoplus_{v \geq e}\Ext^{m}\big(X,\cM,\cN(v)\big).$  We also point out that the integer $r$ is closely related to the Castelnuovo-Mumford regularity of ${_S{N}}$.  This theorem generalize Eisenbud's approximation method for computing sheaf cohomology \cite[theorem~8.3.2]{Loc}.

We prove the main result in two steps.  First, we examine the spectral sequence relating a locally free resolution $\cF_{\bullet}$ of $\cM$ to the global extension modules, namely
$$ E_{2}^{p,q} = \H^{p}\Big( H^{q}\big(X,\SHom_{\cO_{X}}\!(\cF_{\bullet},\cN)\big)\Big) \stackrel{p}{\Longrightarrow} \Ext^{p+q}(X;\cM,\cN).$$
Applying vanishing conditions for sheaf cohomology, we obtain conditions for $r$ which insure that spectral sequence becomes $E_{2}^{m,0} \cong \Ext^{m}(X;\cM,\cN)$.   The problem therefore reduces to computing the cochain complex $\Hom(X;\cF_{\bullet},\cN)$.  In the second step, we assume $\cF_{\bullet}$ corresponds to a free resolution $F_{\bullet}$ of $M_{\geq r}$ and consider the natural homomorphisms from $\Hom_{R}(F_{i},N)$ to $\bigoplus_{v \in \ZZ} \Hom\big(X;\cF_{i},\cN(v)\big)$.  The kernel and cokernel of these maps are both local cohomology modules.  By using local duality, we place further restrictions on $r$ and bound the maximal degree of nonzero elements appearing in the kernel and cokernel.  Combining both parts, we obtain the required isomorphism.

The complexity of the algorithms derived from theorem~\ref{main} depends on computing Gr\"{o}bner bases.  Because of significant differences between worst case bounds for Gr\"{o}bner bases and the complexity of geometric examples, we omit an analysis of the computational complexity.  Instead, we refer the reader to Bayer and Mumford \cite{B-M} for a discussion of the conjectures and results in this area.   

Background material is provided in the first section.  In the second section, we prove two propositions which determine where to truncate $M$.  A proof of theorem~\ref{main} is presented in the third section, followed by algorithms for computing global extension modules and sheaf cohomology.  Section~4 contains four sample calculations: an example in which the bounds in theorem~\ref{main} are sharp, a comparison of two methods for computing the cohomology of a locally free sheaf, an illustration of Serre-Grothendieck duality and the construction of the sheaf associated to an element of $\Ext^{1}(X;\cM,\cN)$.

\section{Preliminaries}

We collect here a number of definitions, standard results and notations.  References for the unproved assertions about free resolutions can be found in Eisenbud \cite{Eis}, the results from algebraic geometry are in Hartshorne \cite{Har}, local cohomology is presented in Brodmann and Sharp \cite{B-S} and the treatment of hypercohomology follows Weibel \cite{Wei}. 

\begin{freeres}
Let $M$ and $N$ be finitely generated graded $R$-modules.  If $e$ is an integer, $M_{e}$ denotes the $e$-th graded component of $M$.  The submodule, $M_{\geq e} = \bigoplus_{d \geq e} M_{d}$, consisting of elements of degree greater than or equal to $e$, is called a truncation of $M$.  For an integer $v$, we write $M(v)$ for the $v$-th twist of $M$ defined by the formula $M(v)_{e} = M_{e+v}$.  When we write isomorphisms and exact sequences of graded modules, we will always arrange that the maps have degree zero.  However, the set of all homogeneous maps from $M$ to $N$ of all degrees is denoted $\Hom_{R}(M,N)$; it is a graded $R$-module, graded by the degrees of the maps.  Using graded free resolutions, this construction extends to a grading on the $R$-modules $\Ext_{R}^{m}(M,N)$.

The (unique) minimal graded free resolution of $M$ will be denoted $F_{\bullet}(M)$ and $\pd_{R}(M)$ will be the length of $F_{\bullet}(M)$, in other words, the projective dimension of $M$.  In particular, there is an exact sequence
\renewcommand{\theequation}{FR}
\begin{equation}
0 \longrightarrow F_{p}(M) \longrightarrow \cdots \longrightarrow F_{1}(M) \longrightarrow F_{0}(M) \longrightarrow 0,
\end{equation}
where $p = \pd_{R}(M)$.  Here each $F_{i}(M)$ is a direct sum of twists of $R$:
$$ F_{i}(M) = \bigoplus_{j=1}^{b_{i}(M)}R\big(-a_{i,j}(M)\big) $$
and hence the maps in (FR) are given by matrices of homogeneous forms.  In this setting, minimality means that none of the entries in these matrices are nonzero constants.

We define $\ol{a}_{i}(M)$ to be the maximal degree of the minimal generators of the $i$-th syzygy module of $M$.  Similarly, $\ul{a}_{i}(M)$ is the minimal degree of the minimal generators of the $i$-th syzygy modules.  Rephrasing, we have
\begin{eqnarray*}
\ol{a}_{i}(M) & = & \max\big\{a_{i,j}(M) : 1 \leq j \leq b_{i}(M) \big\} \\
\ul{a}_{i}(M) & = & \min\big\{a_{i,j}(M) : 1 \leq j \leq b_{i}(M) \big\}.
\end{eqnarray*}
We consider both $\ol{a}_{i}(M)$ and $\ul{a}_{i}(M)$ as functions defined for all integers $i$ by using the conventions $\max\{\emptyset\} = - \infty$ and $\min\{\emptyset\} = \infty$.  We emphasize that $\ol{a}_{i}({_{S}N})$ is determined by a free resolution of ${_{S}N}$ in the category of $S$-modules.
\end{freeres}

\begin{geometry}
By definition, we have $X = \Proj(R)$ and $\PP^{n} = \Proj(S)$.  Sheaves of $\cO_{X}$-modules will usually be denoted by characters in script font, such as $\cG$.  The $\cO_{X}$-module associated to an $R$-module $M$ will be written $\cM$.  We point out that the free resolution $F_{\bullet}(M)$ of the $R$-module $M$ gives rise to a locally free resolution $\cF_{\bullet}(M)$ of the sheaf $\cM$.  The projective embedding of $X$ is given by the morphism of schemes $f \colon X \rightarrow \PP^{n}$.  As $f$ is a closed immersion, the direct image functor $f_{*}$ takes coherent $\cO_{X}$-modules to coherent $\cO_{\PP^{n}}$-modules.  Moreover, $f_{*}\cM$ may be identified with the sheaf obtained from $\cM$ by extending by zero.

The global section functor will be written $\Gamma(X,-)$ and the $q$-th right derived functor, called sheaf cohomology, will be denoted $H^{q}(X,-)$.  The sheaf of local homomorphisms, also called sheaf Hom, is denote $\SHom_{\cO_{X}}\!(\cG,\cG')$ and $\SExt_{\cO_{X}}^{q}\!(\cG,-)$ is the $q$-th right derived functor of $\SHom_{\cO_{X}}\!(\cG,-)$.
\end{geometry}

\begin{local}
The maximal ideal of $S$, which is generated by the elements $x_{0},\ldots,x_{n}$, is called $\nn$.  The image of $\nn$ under the canonical map gives the maximal ideal $\mm$ of $R$.  For each $R$-module $M$, we define $\Gamma_{\mm}(M) = \bigcup_{k \in \NN} \left(0 :_{M} \mm^{k}\right)$; the set of element of $M$ which are annihilated by some power of $\mm$.  Notice that $\Gamma_{\mm}(M)$ is a submodule of $M$.  The $q$-th right derived functor of $\Gamma_{\mm}(-)$ is denoted by $H_{\mm}^{q}(-)$ and is referred to as the $q$-th local cohomology functor with respect to $\mm$.  The local cohomology modules $H_{\mm}^{q}(M)$ are naturally graded Artinian $R$-modules. 
\end{local}

\begin{hypercohomology}
The $p\,$-th right hyper-derived functor of $\Gamma(X,-)$, also called the $p\,$-th hypercohomology functor, is denoted $\HH^{p}(X,-)$.  Let $\cE^{\bullet}$ be a cochain complex of $\cO_{X}$-modules.  We write $\H^{q}(\cE^{\bullet})$ for the $q$-th cohomology group of the complex; that is the kernel of the $q$-th differential modulo the image of the $(q-1)$-th differential.  For a given cochain complex, there are two spectral sequences relating the sheaf cohomology functors to the hypercohomology functors:
\setcounter{equation}{0} 
\renewcommand{\theequation}{SS.\arabic{equation}} 
\begin{eqnarray}
{'\!E}_{2}^{p,q} =  \H^{p} \big( H^{q}(X, \cE^{\bullet} ) \big) & \stackrel{p}{\Longrightarrow} & \HH^{p+q}(X,\cE^{\bullet}) \\
{''\!E}_{2}^{p,q} = H^{p}\big(X, \H^{q}( \cE^{\bullet}) \big) & \stackrel{p}{\Longrightarrow} & \HH^{p+q}(X,\cE^{\bullet}). 
\end{eqnarray}
Both spectral sequences converge when $\cE^{\bullet}$ consists of only a finite number of nonzero $\cO_{X}$-modules.
\end{hypercohomology}

\section{Bounds for Approximating Global Extension Modules}

We begin by giving vanishing conditions for the $m$-th cohomology group of a coherent sheaf.  This lemma can be viewed as refinement of Serre's finiteness theorem \cite[theorem~III.5.2]{Har} or Cartan's theorem~B \cite[section~5.3]{G-H}.

\begin{lemma} \label{vanishing}
Let $N$ be a finitely generated graded $R$-module and let $m$ be a positive integer.  If the inequality $v \geq \ol{a}_{n-m}({_{S}N})-n$ is satisfied, then the cohomology group $H^{m}\big(X,\cN(v)\big)$ vanishes.
\end{lemma}

As hypercohomology will be required in the following discussion, we provide a proof using this machinery.  One could also prove this assertion by induction on the projective dimension of ${_{S}N}$. 

\begin{proof}
Given the closed immersion $f \colon X \rightarrow \PP^{n}$, there exists an isomorphism of cohomology groups \cite[lemma~III.2.10]{Har}:
$$ H^{m}\big(X,\cN(v)\big) \cong H^{m}\big(\PP^{n},f_{*}\cN(v)\big). $$
As $f_{*}\cN(v)$ is the sheaf associated to the $S$-module ${_{S}N}(v)$ \cite[proposition~II.5.12]{Har}, it suffices to prove the assertion for projective space. 

Without loss of generality, we may assume $X = \PP^{n}$, that is $R = S$.  Consider the cochain complex given by $\cE^{-i} = \big(\cF_{i}(N)\big)(v)$.  As $\cF_{\bullet}(N)$ is a resolution of $\cN$, the second hypercohomology spectral sequence (SS.2) collapses and the first spectral sequence (SS.1) becomes
$$ {'\!E}_{2}^{p,q} = \H^{p}\big( H^{q} \big(X,\cE^{\bullet} \big) \big) \stackrel{p}{\Longrightarrow} H^{p+q}\big(X,\cN(v)\big).$$
We recall that, for projective $n$-space, one has the following equivalence \cite[theorem~III.5.1]{Har}:
\renewcommand{\theequation}{VC}
\begin{equation}
H^{q}\big(X,\cO_{X}(v)\big) = 0 \;\; \Longleftrightarrow \;\; \left\{ \renewcommand{\arraycolsep}{1pt}
\begin{array}{lcl}
v < 0 & ; & q = 0 \\
v > -n-1 & ; & q = n \\
v > -\infty & ; & q \neq 0,n 
\end{array} \right] .
\end{equation} 
Because $F_{\bullet}(N)$ is a free resolution of $N$, we have
$$ \cE^{-i} = \bigoplus_{j=1}^{b_{i}(N)} \cO_{X}\big(v-a_{i,j}(N) \big)$$
for $0 \leq i \leq \pd_{R}(N)$; otherwise $\cE^{-i} =0$.  It follows that the relation ${'\!E}_{2}^{p,q} \neq 0$ implies $q$ equals $0$ or $n$ and $-\pd_{R}(N) \leq p \leq 0$.  From this, we see that, for $p+q > 0$, either the target or source of the differential $d_{r}^{p,q} \colon {'\!E}_{r}^{p,q} \rightarrow {'\!E}_{r}^{p+r,q-r+1}$ is zero.  We conclude that $H^{m}\big(X,\cN(v)\big) = \H^{m-n}\big(H^{n}(X,\cE^{\bullet})\big)$.  By making use of the vanishing conditions (VC) for a second time, we observe that, when $-\pd_{R}(N) \leq m-n \leq 0$, the inequality $v \geq \ol{a}_{n-m}(N)-n$ implies that we have $H^{n}\big(X,\cE^{n-m}\big) = 0$.  Otherwise we have $\cE^{n-m} = 0$ and the vanishing is immediate.  Finally, our convention implies that $\ol{a}_{n-m}(N) = - \infty$ for $m-n < - \pd_{R}(N)$ or $m-n > 0$ and the assertion follows.\large{$\Box$}
\end{proof}

\begin{remark} \label{gv} {\rm
We may give a sharper version of the preceding lemma by placing conditions on the integer $m$.  In fact, Grothendieck's vanishing theorem \cite[theorem~III.2.7]{Har} shows that the cohomology group $H^{m}\big(X,\cN(v)\big)$ is zero when $m$ is greater than the dimension of the support of $\cN(v)$, that is $m \geq \dim(N)$ .
} \end{remark}

\begin{example} {\rm
We recall that, for a complete intersection $X$ in $\PP^{n}$ of dimension $d$, one has $H^{m}\big(X,\cO_{X}(v)\big)=0$ for $0 < m < d$ and all integers $v$ \cite[exercise~III.5.5]{Har}.  This illustrates that the converse of Lemma~\ref{vanishing} is false.\large{$\Diamond$} 
} \end{example}

Applying the above lemma, we provide bounds which guarantee that the first hypercohomology spectral sequences (SS.1) collapses.  In these circumstances, we are able to compute $\Ext^{m}\big(X;\cM,\cN(v)\big)$ as the cohomology of a given cochain complex.  More explicitly, we have:

\begin{proposition} \label{collapsing}
Let $M$ and $N$ be two finitely generated graded $R$-modules, let $m$ be an integer and set $\ell = \min\{\dim(N),m\}$.  For any integer $v$ satisfying the inequalities
\begin{eqnarray*} 
v & \geq & \max\big\{ \ol{a}_{n-u}({_{S}N}) -\ul{a}_{m-u}(M) : 1 \leq u \leq \ell \big\}-n, \;\; {\it and} \\
v & \geq & \max\big\{ \ol{a}_{n-u+1}({_{S}N}) - \ul{a}_{m-u}(M) : 2 \leq u \leq \ell \big\}-n,
\end{eqnarray*} 
we have the following isomorphism of graded $R$-modules:
$$ \Ext^{m}\big(X;\cM,\cN(v)\big) \cong \H^{m}\Big( \Hom\big(X;\cF_{\bullet}(M),\cN(v)\big)\Big).$$
\end{proposition}

\begin{proof}
For $m < 0$, both the global extension module and the cohomology group vanish.  Thus, we may confine our attention to $m \geq 0$.  Consider the cochain complex defined by $\cE^{i} = \SHom_{\cO_{X}}\!\big( \cF_{i}(M),\cN(v) \big)$.  Because $\cF_{\bullet}(M)$ is a locally free resolution of $\cM$,  we have \cite[proposition~III.6.5]{Har}:
$$ \H^{q}(\cE^{\bullet}) = \H^{q}\Big(\SHom_{\cO_{X}}\!\big( \cF_{\bullet}(M),\cN(v) \big)\Big) = \SExt_{\cO_{X}}^{q}\!\big(\cM,\cN(v)\big).$$
Furthermore, for any two $\cO_{X}$-modules $\cG$, $\cG'$, Th\'{e}or\`{e}me 4.2.1 in \cite{T} provides the spectral sequence $H^{p}\big(X,\SExt_{\cO_{X}}^{q}\!(\cG,\cG')\big) \stackrel{p}{\Longrightarrow} \Ext^{p+q}(X;\cG,\cG')$.  Hence, the second hypercohomology spectral sequence (SS.2) becomes 
$$ {''\!E}_{2}^{p,q} = H^{p}\Big(X, \SExt_{\cO_{X}}^{q}\!\big(\cM,\cN(v)\big)\Big) \stackrel{p}{\Longrightarrow} \Ext^{p+q}\big(X;\cM,\cN(v)\big) = \HH^{p+q}(X, \cE^{\bullet}),$$  
and consequently the first hypercohomology spectral sequence (SS.1) is
$$ {'\!E}_{2}^{p,q} = \H^{p}\Big(H^{q}\Big(X,\SHom_{\cO_{X}}\!\big(\cF_{\bullet}(M),\cN(v)\big)\Big)\Big) \stackrel{p}{\Longrightarrow} \Ext^{p+q}\big(X;\cM,\cN(v)\big).$$ 
We point out that 
$$H^{0}\Big(X,\SHom_{\cO_{X}}\!\big(\cF_{\bullet}(M),\cN(v)\big)\Big) = \Hom\big(X;\cF_{\bullet}(M),\cN(v)\big),$$ 
whence the relation ${'\!E}_{2}^{p,0} = \H^{p}\Big(\Hom\big(X;\cF_{\bullet}(M),\cN(v)\big)\Big)$.

The proof therefore reduces to showing that the bounds on $v$ imply that the first hypercohomology spectral sequence (SS.1) degenerates; more precisely, we will prove 
$$ \bigoplus_{p+q=m} {'\!E}_{\infty}^{p,q} = {'\!E}_{2}^{m,0}.$$  
To achieve this, it suffices to demonstrate that ${'\!E}_{2}^{m-u,u}$ vanishes for $1 \leq u \leq m$ and ${'\!E}_{2}^{m-u,u-1}$ vanishes for $2 \leq u \leq m$, since the differential $d_{r}^{p,q} \colon E_{r}^{p,q} \rightarrow E_{r}^{p+r,q-r+1}$ has bidegree $(r,1-r)$ and ${'\!E}_{2}^{p,q} = 0$ when $q < 0$ or $p < 0$.  Now, each $\cF_{p}(M)$ is locally free of finite rank, so $\SHom_{\cO_{X}}\!\big( \cF_{p}(M),\cN(v)\big) \cong \bigoplus_{j=1}^{b_{p}(M)} \cN\big(v+a_{p,j}(M)\big)$, and we obtain
$$ {'\!E}_{2}^{p,q} = \H^{p} \left( \bigoplus_{j=1}^{b_{\bullet}(M)} H^{q}\Big(X,\cN\big(v+a_{\bullet,j}(M)\big)\Big)\right).$$
From remark~\ref{gv}, we immediately see that $q \geq \dim(N)$ implies ${'\!E}_{2}^{p,q} = 0$.  Moreover, applying lemma~\ref{vanishing}, we see that if the inequality $v \geq \ol{a}_{n-q}({_{S}N})-\ul{a}_{p}(M) -n$ is satisfied then ${'\!E}_{2}^{p,q}$ is zero.  Combining these observations, the claim follows.\large{$\Box$}
\end{proof}

Next, we relate the $m$-th global extension module to the graded $R$-module $\Ext_{R}^{m}(M,N)$.

\begin{proposition} \label{approx}
Let $M$ and $N$ be two finitely generated graded $R$-modules and let $m$ be an integer.  For any integer $e$ satisfying the inequality
$$ e \geq \max\big\{\ol{a}_{n+1}({_{S}N}),\ol{a}_{n}({_{S}N})\big\} - \ul{a}_{m}(M)-n, $$
there is a natural isomorphism of graded $R$-modules
$$ \Big(\Ext_{R}^{m}(M,N)\Big)_{\geq e} \cong \bigoplus_{v \geq e} \H^{m} \Big( \Hom\big(X;\cF_{\bullet}(M),\cN(v)\big) \Big).$$
\end{proposition}

\begin{proof}
Again, for $m < 0$, both the global extension module and the cohomology group vanish, so we may assume that $m \geq 0$.  Now, the graded $R$-module $\Ext_{R}^{m}(M,N)$ can be computed from a free resolution of $M$, that is 
$$\Ext_{R}^{m}(M,N) \cong \H^{m}\Big( \Hom_{R}\big(F_{\bullet}(M),N\big) \Big).$$  
Relating local cohomology to the sheaf cohomology, we have the following natural exact sequence \cite[20.4.4]{B-S}:
$$ \renewcommand{\arraycolsep}{1pt}
\begin{array}{ccccccc}
0 & \longrightarrow & H_{\mm}^{0} \Big(\Hom_{R}(F_{i}(M),N\big)\Big) & \longrightarrow & \Hom_{R}\big(F_{i}(M),N\big) & \longrightarrow \\
& \longrightarrow & \ds \bigoplus_{v \in \ZZ} \Hom\big(X;\cF_{i}(M),\cN(v) \big) & \longrightarrow & H_{\mm}^{1} \Big(\Hom_{S} \big( F_{i}(M),N \big)\Big) & \longrightarrow & 0,
\end{array}$$
for each $i$; in other words, we have a chain map 
$$ \Hom_{R}\big(F_{\bullet}(M),N\big) \longrightarrow \bigoplus_{v \in \ZZ} \Hom\big(X;\cF_{\bullet}(M), \cN(v)\big).$$
Thus, it suffices to show that this chain map is an isomorphism for degrees greater than or equal to $e$.  

Using the properties of the bifunctor $\Hom_{R}(-,-)$, we have
$$ \Hom_{R}\big(F_{m}(M),N\big) = \Hom_{R}\left( \bigoplus_{j=1}^{b_{m}(M)} R\big(-a_{m,j}(M)\big),N\right) = \bigoplus_{j=1}^{b_{m}(M)} N\big(a_{m,j}(M)\big).$$
Moreover, as $\nn R = \mm$, the $k$-th local cohomology module of the $R$-module $N$, $H_{\mm}^{k}(N)$, is isomorphic to the $k$-th local cohomology module of the $S$-module ${_{S}N}$, $H_{\nn}^{k}({_{S}N})$ \cite[13.1.6]{B-S}.  It therefore suffices to prove that, for $0 \leq k \leq 1$, the maximal degree of an element in $H_{\nn}^{k}({_{S}N})$ is less than $e+\ul{a}_{m}(M)$.  With this in mind, define $\ol{c}_{k}$ to be the maximal degree of an element in $H_{\nn}^{k}({_{S}N})$.  Notice that the local cohomology modules $H_{\nn}^{k}({_{S}N})$ are Artinian, so $\ol{c}_{k}$ is well defined.  Applying the graded form of local duality \cite[13.4.6]{B-S}, we have
$$ H_{\nn}^{k}({_{S}N}) \cong  \Big(\Hom_{\kk}\big( \Ext_{S}^{n+1-k}({_{S}N},S),\kk \big)\Big)(n+1).$$
Setting $\ul{c}_{k}$ to be the minimal degree of an element in $\Ext_{S}^{n+1-k}({_{S}N},S)$, we obtain the formula $\ol{c}_{k} = -\ul{c}_{k}-n-1$.  To compute a lower bound for $\ul{c}_{k}$, we use the free resolution $F_{\bullet}({_{S}N})$ of ${_{S}N}$; in particular, $\Ext_{S}^{k}({_{S}N},S) \cong \H^{k}\Big(\Hom_{S}\big(F_{\bullet}({_{S}N}),S\big) \Big)$. As above, the properties of $\Hom_{S}(-,-)$ give
$$ \Hom_{S}\big( F_{i}({_{S}N}),S \big) \cong \bigoplus_{j=1}^{b_{i}({_{S}N})} S\big( a_{i,j}({_{S}N})\big).$$
Consequently, we have $\ul{c}_{k} \geq -\ol{a}_{n+1-k}({_{S}N})$ which yields the inequality
$$ \max\{\ol{c}_{0},\ol{c}_{1}\} < \max\{\ol{a}_{n}({_{S}N}),\ol{a}_{n+1}({_{S}N})\}-n. $$
The conditions on $e$ clearly imply that $e \geq \max\{\ol{c}_{0},\ol{c}_{1}\} - \ul{a}_{m}(M)$ which completes the proof.\large{$\Box$}
\end{proof}

We end this section by summarizing our results.

\begin{corollary} \label{globalextbound}
Let $M$ and $N$ be two finitely generated graded $R$-modules, let $m$ be an integer and set $\ell = \min\{\dim(N),m\}$.  If the integer $e$ satisfies
 the following inequalities:
\begin{eqnarray*}
e & \geq & \max\big\{ \ol{a}_{n-u}({_{S}N}) -\ul{a}_{m-u}(M) : 0 \leq u \leq \ell \big\}-n; \\
e & \geq & \max\big\{ \ol{a}_{n-u+1}({_{S}N}) - \ul{a}_{m-u}(M) : 0 \leq u \leq \ell, u \neq 1 \big\}-n,  
\end{eqnarray*}
then we have an isomorphism of graded $R$-modules
$$\Big(\Ext_{R}^{m}(M,N)\Big)_{\geq e} \cong \bigoplus_{v \geq e} \Ext^{m}\big(X;\cM,\cN(v)\big).$$
\end{corollary}

\begin{proof}
This follows immediately from proposition~\ref{collapsing} and proposition~\ref{approx}.\large{$\Box$}
\end{proof}

\section{A Method for Computing Global Extension Modules}

Using the bounds developed in the preceding section, we present a proof of the main theorem.

\begin{proofthm}
Fix an integer $r$ such that
$$ r \geq \max\big\{ \ol{a}_{i}({_{S}N})-i : n-\ell \leq i \leq \pd_{S}({_{S}N}) \big\}-m+1.$$
To simplify notation, set $M' = M_{\geq r}$.  We remind the reader that the $\cO_{X}$-module associated to the $R$-module $M$ is equal to the one associated to $M'$ \cite[exercise~II.5.9]{Har}.  Applying corollary~\ref{globalextbound}, it suffices to show that the following inequalities are satisfied:
\renewcommand{\theequation}{I}
\begin{equation}
\left\{ \begin{array}{rclrl}
\ul{a}_{m-u}(M')+n - \ol{a}_{n-u}({_{S}N}) & \geq & 0 & {\rm for} & 0 \leq u \leq \ell;  \\
\ul{a}_{m-u}(M')+n - \ol{a}_{n-u+1}({_{S}N}) & \geq & 0 & {\rm for} & 0 \leq u \leq \ell, u \neq 1.  
\end{array} \right. 
\end{equation}
By convention, $\ul{a}_{m-u}(M') = \infty$ for $m-u < 0$ or $m-u > \pd_{R}(M)$, so we may confine our attention to $m-\pd_{R}(M) \leq u \leq m$.  As $F_{0}(M')$ maps onto $M'$ and the maps $F_{i}(M') \rightarrow F_{i-1}(M')$ are given by matrices of homogeneous forms of positive degree, we have the inequality $\ul{a}_{i}(M') \geq r+i$ for all integers $i$.  By making substitutions for $u$, we transform the inequalities (I) into the following:
$$ \left\{ \begin{array}{rclrcl}
r & \geq & \ol{a}_{i}({_{S}N})-i-m & {\rm for} & n-\ell \leq i \leq n;  \\
r & \geq & \ol{a}_{i}({_{S}N})-i-m+1 & {\rm for} & n-\ell+1 \leq i \leq n+1,& i \neq n.
\end{array} \right. $$
Again, our conventions imply that $\ol{a}_{i}({_{S}N}) = -\infty$ for $i < 0$ or $i > \pd_{S}({_{S}N})$, so we need only consider $0 \leq i \leq \pd_{S}({_{S}N})$.  The choice of $r$ therefore guarantees that these inequalities hold.\large{$\Box$}
\end{proofthm}

\begin{remark} \label{dim0} {\rm
For $\dim(M) =0$, it is necessary and sufficient that $M$ have finite length \cite[corollary~2.17]{Eis}.  In particular, if $\dim(M)=0$ then $M$ has only finitely many nonzero graded components and we have $M_{\geq r} =$ for $r \gg 0$.  It follows from theorem~\ref{main} that $\Ext^{m}\big(X;\cM,\cN\big)=0$ when $\dim(M) =0$.  
} \end{remark}

For the special case $X = \PP^{n}$, we have the following additional vanishing result:

\begin{corollary}
Let $\cM$ and $\cN$ be two coherent $\cO_{X}$-modules and assume $X = \PP^{n}$. If $m > n$, then the global extension module $\Ext^{m}(X;\cM,\cN)$ vanishes.
\end{corollary}

\begin{proof}
Let $R$ be the homogeneous coordinate ring for $X$, and let $M$, $N$ be two finitely generated graded $R$-modules which correspond to $\cM$, $\cN$ respectively.  By theorem~\ref{main}, we are required to show that, for $m > n$, the module $\big(\Ext_{R}^{m}(M_{\geq r},N)\big)_{0}$ vanishes for $r$ sufficiently large.  It therefore suffices to prove that $\pd_{R}(M_{\geq r}) \leq n$ for $r \gg 0$.  Applying remark~\ref{dim0}, we may assume $\dim(M) > 0$.  Now, since $X= \PP^{n}$, we have $R = S = \kk[x_{0},\ldots,x_{n}]$ and the Hilbert syzygy theorem \cite[theorem~1.13]{Eis} provides the inequality $\pd_{R}(M_{\geq r}) \leq n+1$.  Moreover, as $R$ is a Cohen-Macaulay ring, we have $\pd_{R}(M_{\geq r}) = \dim(R) = n+1$ if and only if the maximal ideal $\mm$ of $R$ is associated to $M_{\geq r}$ \cite[corollary~19.10]{Eis}.  With this in mind, recall that $H_{\mm}^{0}(M)$ is the set of elements of $M$ which are annihilated by some power of $\mm$.  Since the module $H_{\mm}^{0}(M)$ is Artinian, there exists a maximal degree $\ol{c}_{0}$ of an element in $H_{\mm}^{0}(M)$.  Hence, if $r > \ol{c}_{0}$ then the maximal ideal $\mm$ is not associated to $M_{\geq r}$ which completes the argument.\large{$\Box$}
\end{proof}

Using Gr\"{o}bner bases, one can compute the minimal free resolution of a finitely generated graded $S$-module (see chapter~15 in Eisenbud \cite{Eis}).  This allows one to calculate the projective dimension of a finitely generated graded $S$-module and the degrees of the minimal generators of the syzygy modules.  Furthermore, using the free resolutions one can construct the Ext modules for a pair of finitely generated graded $R$-modules --- again see chapter~15 in Eisenbud \cite{Eis}.  Grayson and Stillman \cite{M2} have implemented these functions in their software system, {\tt Macaulay2}.  With this in mind, theorem~\ref{main} yields algorithm~\ref{globalExtSum}.  By extracting the appropriate graded part, we have algorithm~\ref{globalExt} which returns the $m$-th global extension module.  For the special case  $M=R$, we have the equality \cite[proposition~III.6.4]{Har}:
\renewcommand{\theequation}{EH}
\begin{equation}
\Ext^{m}\big(X;\cO_{X},\cN(v)\big) \cong H^{m}\big(X,\cN(v)\big).
\end{equation}
Hence, we obtain the following two techniques for computing sheaf cohomology: algorithm~\ref{sheafCohomologySum} and algorithm~\ref{sheafCohomology}.  These algorithms have been implemented in {\tt Macaulay2}; the routines may be found at {\tt http://math.berkeley.edu/$\widetilde{\;\;}$ggsmith}. 

\begin{figure}
\fbox{
\begin{minipage}[c]{28pc} 
\begin{algorithm} \label{globalExtSum} ${\tt globalExtSum\/}(m,e,M,N)$ \hfill
{\footnotesize 
\begin{tabbing}
{\bf input}: \= \kill 
{\bf input}: \> two integers $m$, $e$ and two finitely generated graded $R$-modules $M$, $N$. \\ 
{\bf assumption}: \= \kill
{\bf assumption}: \> the graded ring $R$ is a quotient of a polynomial ring $S$. \\
{\bf output}: \= \kill 
{\bf output}: \> the graded $R$-module $\bigoplus_{v \geq e} \Ext^{m}\big(X;\cM,\cN(v)\big)$, where $X = \Proj(R)$ and \\
\> $\cM$, $\cN$ are the coherent $\cO_{X}$-modules associated to $M$, $N$ respectively. \\
{\bf begin} \\ 
iiiii \= iiiii \= iiiii \= iiiii \= \kill 
\> {\bf if} $\dim(M) = 0$ {\bf or} $m < 0$ {\bf then} $E : = 0$ \\
\> {\bf else} $($ \\
\> \> $n :=$ the number of generators of $S$ minus $1$; \\
\> \> $\ell := \min\{\dim(N),m\}$; \\
\> \> ${_{S}N} :=$ the $S$-module obtained from $N$ by restriction of scalars; \\
\> \> $F_{\bullet}({_{S}N}) :=$ the minimal free resolution of ${_{S}N}$; \\
\> \> $\pd({_{S}N}) :=$ the length of $F_{\bullet}({_{S}N})$; \\
\> \> {\bf if} $\pd({_{S}N}) < n-\ell$ {\bf then} $E := \Ext_{R}^{m}(M,N)$ \\
\> \> {\bf else} $($ \\
\> \> \> {\bf for} $i$ {\bf from} $n-\ell$ {\bf to} $\pd({_{S}N})$ {\bf do} $($ \\ 
\> \> \> \> $\ol{a}_{i}(N) :=$ the maximal degree of the minimal generators of $F_{i}({_{S}N})$; $);$ \\
\> \> \> $r := \max\big\{ \ol{a}_{i}({_{S}N}) : n-\ell \leq i \leq \pd({_{S}N}) \big\}-e-m+1$; \\
\> \> \> $E := \Ext_{R}^{m}(M_{\geq r},N)$; $);$ $);$ \\
\> $e' :=$ minimal degree of the generators of $E$; \\
\> {\bf if} $e' < e$ {\bf then} $E := E_{\geq e}$; \\
\> {\bf return} $E$; \\
{\bf end}.  
\end{tabbing}}
\end{algorithm}
\end{minipage} }
\end{figure}

The complexity of algorithm~\ref{globalExtSum} will depend on computing $F_{\bullet}({_{S}N})$ and $F_{\bullet}(M_{\geq r})$; the second resolution is used to determine $\Ext_{R}^{m}(M_{\geq r},N)$.  To find these resolutions involves several Gr\"{o}bner bases calculations.  As noted in the introduction, there is a significant difference between worst case bounds for Gr\"{o}bner bases and the complexity of geometric examples.  With this in mind, we omit an analysis of complexity and refer the reader to Bayer and Mumford \cite{B-M} for a discussion of the conjectures and results in this area.

\begin{figure}
\fbox{
\begin{minipage}[c]{28pc} 
\begin{algorithm} \label{globalExt}  ${\tt globalExt\/}(m,M,N)$ \hfill
{\footnotesize 
\begin{tabbing}
{\bf input}: \= \kill 
{\bf input}: \> an integer $m$ and two finitely generated graded $R$-modules $M$, $N$. \\
{\bf assumption}: \= \kill
{\bf assumption}: \> the graded ring $R$ is a quotient of a polynomial ring $S$. \\
{\bf output}: \= \kill 
{\bf output}: \> the $\kk$-vector space $\Ext^{m}(X;\cM,\cN)$, where $X = \Proj(R)$ and $\cM$, $\cN$ are the \\
\> coherent $\cO_{X}$-modules associated to $M$, $N$ respectively. \\
{\bf begin} \\ 
iiiii \= iiiii \= iiiii \= iiiii \= $\ol{a}_{i}(N) :=$ \= \kill 
\> $E := {\tt globalExtSum\/}(m,0,M,N)$; \\
\> $V :=$ the $\kk$ vector space isomorphic to $E_{0}$;  \\
\> {\bf return} $V$; \\
{\bf end}.  
\end{tabbing}}
\end{algorithm}
\end{minipage} }
\end{figure}

\begin{figure}
\fbox{
\begin{minipage}[c]{28pc} 
\begin{algorithm} \label{sheafCohomologySum} ${\tt sheafCohomologySum\/}(m,e,N)$ \hfill
{\footnotesize 
\begin{tabbing}
{\bf input}: \= \kill 
{\bf input}: \> two integers $m$, $e$ and a finitely generated graded $R$-module $N$. \\
{\bf assumption}: \= \kill
{\bf assumption}: \> the graded ring $R$ is a quotient of a polynomial ring $S$. \\
{\bf output}: \= \kill 
{\bf output}: \> the graded $R$-module $\bigoplus_{v \geq e} H^{m}\big(X,\cN(v)\big)$, where $X = \Proj(R)$ and $\cN$ is the \\
\> coherent $\cO_{X}$-module associated to $N$. \\
{\bf begin} \\ 
iiiii \= iiiii \= iiiii \= iiiii \= $\ol{a}_{i}(N) :=$ \= \kill 
\> $E := {\tt globalExtSum\/}(m,e,R,N)$; \\
\> {\bf return} $E$; \\
{\bf end}.  
\end{tabbing}}
\end{algorithm}
\end{minipage} }
\end{figure}

\begin{figure}
\fbox{
\begin{minipage}[c]{28pc} 
\begin{algorithm} \label{sheafCohomology}  ${\tt sheafCohomology\/}(m,N)$ \hfill
{\footnotesize 
\begin{tabbing}
{\bf input}: \= \kill 
{\bf input}: \> an integer $m$ and a finitely generated graded $R$-module $N$. \\
{\bf assumption}: \= \kill
{\bf assumption}: \> the graded ring $R$ is a quotient of a polynomial ring $S$. \\
{\bf output}: \= \kill 
{\bf output}: \> the $\kk$-vector space $H^{m}\big(X,\cN(v)\big)$, where $X = \Proj(R)$ and $\cN$ is the coherent \\
\> $\cO_{X}$-module associated to $N$. \\
iiiii \= \kill 
{\bf begin} \\ 
\> $E := {\tt globalExt\/}(m,R,N)$; \\
\> {\bf return} $E$; \\
{\bf end}.  
\end{tabbing} }
\end{algorithm}
\end{minipage} }
\end{figure}

\begin{remark} {\rm 
Algorithm~\ref{sheafCohomologySum} is equivalent to the approximation method for computing sheaf cohomology described in theorem C.3.1 in \cite{Um2}.
} \end{remark}

\section{Some Examples}

All computations in this section were done on an Intel Pentium 166~MHz machine running {\tt Macaulay2} (version 0.8.41) for the Linux platform.

We begin by giving an example in which the bounds in theorem~\ref{main} are sharp. 

\begin{example} {\rm
Let $X = \PP^{3}$, let $\cM = \cO_{X}$, let $\cN$ be the structure sheaf of the smooth rational quartic curve in $\PP^{3}$ and set $m=1$.  Following Eisenbud, Grayson and Stillman \cite{Um2}, these construction in {\tt Macaulay2} are:
{\footnotesize \begin{verbatim}
    i1 : kk = ZZ/32003; S = kk[w,x,y,z];

    i3 : I = monomialCurve(S,{1,3,4})

                          3    2     2  2    3  2
    o3 = ideal (x*y-w*z, y -x*z , w*y -x z, x -w y)

    o3 : Ideal of S

    i4 : N = S^1/I

    o4 = cokernel {0} | xy-wz y3-xz2 wy2-x2z x3-w2y |

                                  1
    o4 : S - module, quotient of S
\end{verbatim}} 
\noindent Next, we determine $\max\big\{ \ol{a}_{i}({_{S}N})-i : n-\ell \leq i \leq \pd_{S}({_{S}N}) \big\}$.  Recall the minimal free resolution of $N$ \cite[8.1]{Loc}:
$$ 0 \longrightarrow S(-5) \longrightarrow \bigoplus_{j=1}^{4} S(-4) \longrightarrow S(-2) \oplus \left( \bigoplus_{j=1}^{2} S(-3) \right) \longrightarrow S;$$
as $n=3$ and $\ell = 1$, we obtain $r \geq 2$.  We verify that {\tt globalExtSum((1,0),S,N)} agrees with $\big(\Ext_{S}^{1}(S_{\geq 2},N)\big)_{\geq 0}$.     
{\footnotesize \begin{verbatim}
    i5 : globalExtSum((1,0),S,N)  

    o5 = 0

    o5 : S - module

    i6 : prune truncate(0, Ext^1(truncate(2,S^1),N)) 

    o6 = 0

    o6 : S - module
\end{verbatim}}
\noindent Lastly, we show that $\big(\Ext_{S}^{1}(S_{\geq 1},N)\big)_{\geq 0} \neq 0$.
{\footnotesize \begin{verbatim}
    i7 : prune truncate(0, Ext^1(truncate(1,S^1),N)) 

    o7 = cokernel {0} | 0 0 0 z 0 y 0 0 0 0 x 0 0 0 w 0 |
                  {0} | z 0 0 0 0 0 y 0 0 0 0 x 0 0 0 w |
                  {0} | 0 z 0 0 y 0 0 0 x 0 0 0 w 0 0 0 |
                  {0} | 0 0 z 0 0 0 0 y 0 x 0 0 0 w 0 0 |

                                  4
    o7 : S - module, quotient of S
\end{verbatim}}
\noindent We conclude that the restrictions on $r$ in theorem~\ref{main} are optimal.\large{$\Diamond$}
} \end{example}

\begin{cohomology}
Recall that, for a locally free sheaf $\cL$ on $X$, there is an isomorphism \cite[propostion~III.6.7]{Har}:
$$ \Ext^{m}(X;\cL^{\vee},\cO_{X}) \cong \Ext^{m}(X;\cO_{X},\cL),$$
where $\cL^{\vee} = \SHom_{\cO_{X}}\!(\cL,\cO_{X})$.  It follows from the equation (EH) that both of these global extension modules are isomorphic to $H^{m}(X,\cL)$.  In particular, this gives two different methods for computing the cohomology of a locally free sheaf.  In the next example, we compute the cohomology of a specific locally free sheaf by these two different approaches and compare the CPU time used.
\end{cohomology}

\begin{example} {\rm
Consider the Veronese surface $X$ embedded in $\PP^{5}$.  Following example~2.6 in Harris \cite{harris}, $X$ is defined by the $2 \times 2$ minors of a generic $3 \times 3$ symmetric matrix.  In {\tt Macaulay2}, we construct the homogeneous coordinate ring $R$ of $X$ as follows: 
{\footnotesize \begin{verbatim}
    i1 : kk = ZZ/32003; 

    i2 : S = kk[u,v,w,x,y,z];

    i3 : I = minors(2,genericSymmetricMatrix(S,u,3))

                  2                                      2
    o3 = ideal (-v +u*x, -v*w+u*y, -w*x+v*y, -v*w+u*y, -w +u*z,
                                                                  2
                                  -w*y+v*z, -w*x+v*y, -w*y+v*z, -y +x*z)

    o3 : Ideal of S

    i4 : R = S/I;
\end{verbatim}}
\noindent Because $X$ is smooth \cite[exercise~14.13]{harris}, the cotangent bundle $\Omega$ is a locally free sheaf \cite[theorem~II.8.15]{Har}.  Moreover, Eisenbud, Grayson and Stillman \cite{Um2} provide an algorithm called {\tt cotangentBundle} which returns the module representing the cotangent bundle.  Using this, we create $\Omega$ and $\Omega^{\vee}$:
{\footnotesize \begin{verbatim}
    i5 : Omega = cotangentBundle R;
       
    i6 : OmegaDual = dual Omega;
\end{verbatim}}
\noindent We establish the identity $\bigoplus_{v \geq 0} \Ext^{1}\big(X;\Omega^{\vee},\cO_{X}(v)\big) \cong \bigoplus_{v \geq 0} \Ext^{1}\big(X;\cO_{X},\Omega(v)\big)$.  Evaluating the left hand side, we have:
{\footnotesize \begin{verbatim}
    i7 : LHS = globalExtSum((1,0),OmegaDual,R)
 
    o7 = cokernel {0} | z y x w v u |

                                  1
    o7 : R - module, quotient of R
\end{verbatim}}
\noindent The {\tt Macaulay2} output for the right hand side takes the form:
{\footnotesize \begin{verbatim}
    i8 : RHS = globalExtSum((1,0),R,Omega)
 
    o8 = cokernel {0} | z y x w v u |

                                  1
    o8 : R - module, quotient of R
\end{verbatim}}
\noindent Next, we use the {\tt benchmark} function to produce an accurate timing for these techniques.  The output is the number of seconds used by the CPU.
{\footnotesize \begin{verbatim}
    i9 : benchmark("globalExtSum((1,0),OmegaDual,R)")

    o9 = 0.51

    o9 : RR

    i10 : benchmark("globalExtSum((1,0),R,Omega)")

    o10 = 516.18

    o10 : RR
\end{verbatim}}
\noindent One immediately sees the significant difference in the length of two calculations.  For this particular example, we can explain the difference by examining where the first module is truncated.  For the left hand side, no truncation is required since $\pd_{S}({_{S}R})=3<4=n-\ell$.  Hence, evaluating {\tt globalExtSum((1,0),OmegaDual,R)} involves finding the minimal free resolution of the $R$-module $\Omega$ and the $S$-module ${_{S}R}$.  However, for the right hand side, we truncate at $r=2$; thus calculating {\tt globalExtSum((1,0),R,Omega)} requires the minimal free resolution of $R$-module $R_{\geq 2}$ and the $S$-module ${_{S}\Omega}$.  Now, the number of minimal generators of $R_{\geq 2}$ is significantly larger than the number for ${_{S}R}$ and this leads to the difference in the amount of time used.  From this example, we conclude that, at least for locally free sheaves, algorithm~\ref{sheafCohomologySum} and algorithm~\ref{sheafCohomology} are not the most efficient methods.\large{$\Diamond$}
} \end{example}

\begin{duality} 
For the third application, we illustrate Serre-Grothendieck duality.  Recall that if $X$ is Cohen-Macaulay, closed subscheme of $\PP^{n}$ of pure dimension $d$, then there exists isomorphisms $\Ext^{d-j}\big(X;\cG,\omega_{X}) \cong H^{j}(X,\cG)^{\vee}$ for $0 \leq j \leq d$, where $\cG$ is a coherent sheaf on $X$, $\omega_{X}$ is the dualizing sheaf and $^{\vee}$ denotes the dual vector space \cite[theorem~III.7.6]{Har}.
\end{duality}

\begin{example} {\rm
Let $X$ be the Del Pezzo surface of degree $4$ in $\PP^{4}$; $X$ is a complete intersection of two quadratic hypersurfaces in $\PP^{4}$ \cite[exercise~III.4.13]{Har}.  We define the homogeneous coordinate ring $R$ of $X$ as follows:
{\footnotesize \begin{verbatim}
    i1 : kk = ZZ/32003; 

    i2 : S = kk[v,w,x,y,z];

    i3 : I = ideal(w*x,y*z);

    o3 : Ideal of S
       
    i4 : R = S/I;
\end{verbatim}}
\noindent By definition, we know that the dualizing sheaf $\omega_{X}$ is isomorphic to $\cO_{X}(-1)$ \cite[remark~III.4.7]{Har}.  For simplicity, we let $\cG$ be the sheaf corresponding to the cokernel of a generic symmetric matrix.  In {\tt Macaulay2} this appears as:
{\footnotesize \begin{verbatim}
    i5 : omega = R^{-1};

    i6 : G = coker genericSymmetricMatrix(R,v,2) 

    o6 = cokernel {0} | v w |
                  {0} | w x |

                                  2
    o6 : R - module, quotient of R
\end{verbatim}}
\noindent Finally, we evaluate both sides of the duality equation.  For $j=0$, we have
{\footnotesize \begin{verbatim}
    i7 : EE2 = globalExt(2,G,omega)
     
           2
    o7 = kk

    o7 : kk - module, free

    i8 : HH0 = Hom(sheafCohomology(0,G),kk)
     
           2
    o8 = kk

    o8 : kk - module, free
\end{verbatim}}
\noindent In the case $j=1$, we obtain 
{\footnotesize \begin{verbatim}
    i9 : EE1 = globalExt(1,G,omega)
      
           2
    o9 = kk

    o9 : kk - module, free

    i10 : HH1 = Hom(sheafCohomology(1,G),kk)
       
            2
    o10 = kk

    o10 : kk - module, free
\end{verbatim}}
\noindent and, when $j = 0$, it is
{\footnotesize \begin{verbatim}
    i11 : EE0 = globalExt(0,G,omega)
        
    o11 = 0

    o11 : kk - module

    i12 : HH2 = Hom(sheafCohomology(2,G),kk)
        
    o12 = 0
\end{verbatim}}
\noindent This completes the example.\large{$\Diamond$}
} \end{example}

\begin{extensions} 
As a final example, we compute a nontrivial global extension.  We reminder the reader that a global extension of $\cM$ by $\cN$ is an exact sequence of $\cO_{X}$-modules:  $0 \longrightarrow \cN \longrightarrow \cE \longrightarrow \cM \longrightarrow 0$;  two extension are equivalent if there is a commutative diagram
$$ \renewcommand{\arraycolsep}{1pt}
\begin{array}{ccccccccc}
0 & \longrightarrow & \cN & \longrightarrow & \cE & \longrightarrow & \cM & \longrightarrow & 0 \\
& & \big\downarrow\vcenter{\rlap{${\scriptstyle id_{\cN}}$}} & & \big\downarrow\vcenter{\rlap{${\scriptstyle \cong}$}} & & \big\downarrow\vcenter{\rlap{${\scriptstyle id_{\cM}}$}} \\
0 & \longrightarrow & \cN & \longrightarrow & \cE' & \longrightarrow & \cM & \longrightarrow & 0.
\end{array} $$
The global extension modules derive their name from the following:  the equivalence classes of global extensions are in bijective correspondence with $\Ext^{1}(X;\cM,\cN)$ \cite[section~5.3]{G-H}.

Applying theorem~\ref{main}, there is a bijection between equivalence classes of global extensions and $\big(\Ext_{R}^{1}(M_{\geq r},N)\big)_{0}$ for $r$ sufficiently large.  Moreover, the functor which takes a module to its associated sheaf is exact, so each global extension 
$$0 \longrightarrow \cN \longrightarrow \cE \longrightarrow \cM \longrightarrow 0$$ 
corresponds to an exact sequence of $R$-modules 
$$0 \longrightarrow N \longrightarrow E \longrightarrow M_{\geq r} \longrightarrow 0,$$ 
again for $r \gg 0$.
\end{extensions}

\begin{example} {\rm  
Let $X$ be the plane elliptic curve defined by $x^3+y^3=z^3$.  We will construct the unique nontrivial extension of $\cO_{X}$ by $\cO_{X}$.  We divide this computation into seven parts.
\begin{enumerate}
\item[(i)]  We create the homogeneous coordinate ring $R$ of $X$.
{\footnotesize \begin{verbatim}
    i1 : kk = ZZ/32003; 

    i2 : S = kk[x,y,z];

    i3 : I = ideal (x^3+y^3-z^3);

    o3 : Ideal of S

    i4 : R = S/I

    o4 = R

    o4 : QuotientRing
\end{verbatim}}
\item[(ii)]  We calculate the $\kk$-vector space $\Ext^{1}(X,\cO_{X},\cO_{X})$:
{\footnotesize \begin{verbatim}
    i5 : globalExt(1,R,R)

           1
    o5 = kk

    o5 : kk - module, free
\end{verbatim}}
\noindent As $\Ext^{1}(X;\cO_{X},\cO_{X})$ is one dimensional, there exists only one nontrivial global extension up to equivalence.  
\item[(iii)]  We determine where to truncate the first module.  Following algorithm~\ref{globalExtSum}, we first calculate $n$ and $\ell$.
{\footnotesize \begin{verbatim}
    i6 : n = (numgens S)-1, l = min(dim(R),1)

    o6 = (2, 1)

    o6 : Sequence
\end{verbatim}}
\noindent We see that one must choose $r \geq \max\big\{ \ol{a}_{i}({_{S}N})-i : 1 \leq i \leq \pd_{S}({_{S}N})\big\}$.  After computing $F_{\bullet}({_{S}N})$, we simply set $r$ equal to the maximum: 
{\footnotesize \begin{verbatim}
    i7 : sN = S^1/I

    o7 = cokernel {0} | x3+y3-z3 |

                                  1
    o7 : S - module, quotient of S

    i8 : FF = res sN
       
          1      1
    o8 = S  <-- S
             
         0      1

    o8 : ChainComplex

    i9 : r = (max degrees FF_1)#0-1

    o9 = 2
\end{verbatim}}
\noindent We introduce the modules $M'$ and $N$ to help avoid confusing notation.
{\footnotesize \begin{verbatim}
    i10 : M' = prune truncate(2,R^1)

    o10 = cokernel {2} | 0  y  0  -z 0  0  0  0  x  |
                   {2} | y  -z 0  0  0  0  0  x  0  |
                   {2} | -z 0  0  y  0  x  0  0  0  |
                   {2} | 0  0  y  0  0  0  x  0  -z |
                   {2} | 0  0  -z 0  x  -y 0  -z 0  |
                   {2} | 0  0  0  x  -y 0  -z 0  0  |

                                   6
    o10 : R - module, quotient of R

    i11 : N = R^1

           1
    o11 = R

    o11 : R - module, free
\end{verbatim}} 
\item[(iv)]  We fix an exact sequence $0 \longrightarrow K \stackrel{\alpha}{\longrightarrow} P \stackrel{\beta}{\longrightarrow} M' \longrightarrow 0$, where $P$ a free $R$-module.  This sequence is constructed from a free presentation of the $R$-module $M'$, that is $F_{1} \stackrel{\mu}{\longrightarrow} F_{0} \longrightarrow M' \longrightarrow 0$.
{\footnotesize \begin{verbatim}
    i12 : mu = presentation M'; P = target mu; beta = map(M', P, id_P);

    i14 : K = image mu; alpha = map(P,K,mu);
\end{verbatim}}
\noindent We verify that we have an short exact sequence.
{\footnotesize \begin{verbatim}
    i17 : {ker alpha == 0, image alpha == ker beta, image beta == M'}

    o17 = {true, true, true}

    o17 : List
\end{verbatim}}
\noindent Now, this short exact sequence yields an exact sequence
$$ \Hom_{R}(P,N) \longrightarrow \Hom_{R}(K,N) \longrightarrow \Ext_{R}^{1}(M',N) \longrightarrow 0.$$
Thus, an element in $\Ext_{R}^{1}(M,N)$ gives rise to an element $\theta$ in $\Hom_{R}(K,N)$.
\item[(v)]  We choose a homomorphism $\theta \colon K \rightarrow N$ directly.  
{\footnotesize \begin{verbatim}
    i18 : morphisms = Hom(K, N)

    o18 = image {-3} | 0  y  0  -z 0  0  0   x2  0   |
                {-3} | y  -z 0  0  0  0  x2  0   xz  |
                {-3} | 0  0  y  0  -z 0  0   0   x2  |
                {-3} | -z 0  0  y  0  x  0   0   0   |
                {-3} | 0  0  0  0  x  -y 0   0   -z2 |
                {-3} | 0  0  0  x  -y 0  -yz -z2 0   |
                {-3} | 0  0  x  0  0  -z 0   0   -y2 |
                {-3} | 0  x  0  0  -z 0  -z2 -y2 0   |
                {-3} | x  0  -z 0  0  0  -y2 -yz 0   |

                                    9
    o18 : R - module, submodule of R

    i19 : choice = matrix{{1_R},{0},{0},{0},{0},{0},{0},{0},{0}};
 
    i20 : theta = homomorphism map(morphisms, R^1, choice)

    o20 =  | 0 y 0 -z 0 0 0 0 x |

    o20 : Matrix
\end{verbatim}} 
\noindent The module $E$ is then the pushout of $\theta$ and $\alpha$; more explicitly the cokernel of the map $\psi \colon K \rightarrow P \oplus N$ where $\psi(x) = \big(\alpha(x),-\theta(x)\big)$.
\item[(vi)]  We construct the module $E$.
{\footnotesize \begin{verbatim}
    i21 : D = P ++ N;

    i22 : psi = map(D, K, alpha || -theta)

    o22 =  | 0  y  0  -z 0  0  0  0  x  |
           | y  -z 0  0  0  0  0  x  0  |
           | -z 0  0  y  0  x  0  0  0  |
           | 0  0  y  0  0  0  x  0  -z |
           | 0  0  -z 0  x  -y 0  -z 0  |
           | 0  0  0  x  -y 0  -z 0  0  |
           | 0  -y 0  z  0  0  0  0  -x |

    o22 : Matrix

    i23 : E = coker psi;
\end{verbatim}}
\item[(vii)]  We check that we have the short exact sequence $0 \longrightarrow N \stackrel{\iota}{\longrightarrow} E \stackrel{\phi}{\longrightarrow} M' \longrightarrow 0$.
{\footnotesize \begin{verbatim}
    i24 : iota = map(E, N, map(P,N,0) || id_N);  

    i25 : phi = map(M', E, id_P | map(P,N,0));

    i26 : {ker iota == 0, image iota == ker phi, image phi == M'}

    o26 = {true, true, true}

    o26 : List
\end{verbatim}}
\end{enumerate}
\noindent Therefore, the $R$-module $E$ represents the unique nontrivial global extension.  For completeness, we calculate the global sections of the sheaf associated to $E$.
{\footnotesize \begin{verbatim}
    i27 : sheafCohomology(0,E)

            1
    o27 = kk

    o27 : kk - module, free
\end{verbatim}}
\noindent As the space of global sections has dimension one, it follows that $E$ represents an indecomposable rank 2 vector bundle \cite[theorem~III.2.15]{Har}.\large{$\Diamond$}
} \end{example}

\section*{Acknowledgements}

This work started as a project for Richard Fateman's course on algebraic algorithms.  I would like to thank my advisor David Eisenbud for suggesting this problem and for many useful discussions. 

%  ****************************************
%  *         BIBLIOGRAPHY                 *
%  ****************************************
%

\end{document}